\documentclass[11pt]{article}

\usepackage[english]{babel}
\usepackage{amsfonts,amssymb,amsmath,amsthm,euscript,xspace,cite}
\usepackage{graphicx}
\usepackage{xcolor}
\usepackage{secdot}
\sectiondot{subsection}

\theoremstyle{plain}

\newcommand{\hmove}[1]{\overset{#1}{\longrightarrow}}
\newtheorem{theorem}{Theorem}

\newtheorem{lemma}[theorem]{Lemma}
\newtheorem{proposition}[theorem]{Proposition}
\newtheorem{corollary}[theorem]{Corollary}
\theoremstyle{definition}
 
\newtheorem{remark}{Remark}

\def\ZZ{\mathbb Z}

\def\G{\ensuremath{\mathcal G}}
\def\cG{\ensuremath{\mathcal G}}

\def\R{\ensuremath{\mathcal R}}

\def\cH{\mathcal{H}}

\def\mim#1{\ensuremath{\mathrm {NIM}^{\leq}_{#1}}}

\title{Computing Remoteness Functions 
of Moore, Wythoff, and Euclid's games}
\begin{document}
\maketitle
\begin{center}
\uppercase{\bf }
\vskip 20pt

{\bf Endre Boros}\\
{MSIS and RUTCOR, RBS, Rutgers University, 
NJ,  USA}\\
{\tt endre.boros@rutgers.edu}\\ 
\vskip 7pt
{\bf Vladimir Gurvich}\\
{RUTCOR, Rutgers University, NJ, USA;  \\
National Research University Higher School of Economics,  Moscow, Russia}\\
{\tt vladimir.gurvich@gmail.com}\\ 
\vskip 7pt
{\bf Kazuhisa Makino}\\
{Research Institute for Mathematical Sciences, 
Kyoto University, 
Kyoto, 
Japan}\\
{\tt makino@kurims.kyoto-u.ac.jp}\\ 
\vskip 7pt
{\bf Michael Vyalyi}\\
{National Research University Higher School of Economics, Moscow, Russia;\\
Institute of Physics and Technology, Dolgorpudnyi, Russia;\\
Federal Research Center ``Computer Science and Control" 
of the Russian Academy of Sciences, Moscow, Russia}\\
{\tt vyalyi@gmail.com}
\end{center}

\maketitle

\begin{abstract}
We study remoteness function $\R$ of impartial games 
introduced by Smith in 1966. 
The player who moves from a position $x$  
can win if and only if  $\R(x)$  is odd.  
The odd values of $\R(x)$ show how soon the winner can win,
while even values show how long the loser can resist,
provided both players play optimally. 
This function can be applied to 
the conjunctive compounds of impartial games, 
in the same way as the Sprague-Grundy function 
is applicable to their disjunctive compounds. 
\newline 
We provide polynomial algorithms computing 
$\R(x)$  for games Euclid and generalized Wythoff. 
For Moore's NIM we give a simple explicit formula 
for $\R(x)$  if it is even and show that 
computing it  becomes an NP-hard problem for the odd values.
\end{abstract}

\section{Introduction}\label{s1}
\label{s0}
\subsection{Impartial Games}
We assume that the reader is familiar 
with basic concepts of impartial game theory; 
see e.g., \cite{ANW07,BCG01-04} for an introduction.  

An impartial game is modeled by an acyclic directed graph 
(digraph)  $\Gamma$  whose vertices and directed edges 
are called, respectively, 
{\em positions} and {\em moves} of the corresponding game.
We assume that  $\Gamma$  is  
{\em potentially finite}, that is, 
for any fixed initial position  $x^0$ 
the set of positions that can be reached from $x^0$ 
by a sequence of moves is finite. 
A position with no moves is called {\em terminal}. 
The player who has to move in it 
(but cannot) loses the {\em normal} version 
and wins its {\em mis\`ere} version. 
In this paper, we restrict ourselves by the normal version. 

\subsection{Smith's Remoteness Functions} 
\label{ss00}
In 1966 Smith \cite{Smi66} 
introduced a {\em remoteness} function $\R$ 
for impartial games by the following algorithm. 

\smallskip 

Set  $\R = 0$  for all terminal positions 
of  $\Gamma$  and  
$\R(x) = 1$  if and only if 
there is a move from $x$ to a terminal position.  
Delete all labeled positions from  $G$  and 
repeat the above procedure increasing  $\R$  by 2,  
that is, assign 2 and 3  instead of  0 and 1; etc. 

This algorithm was considered 
as early as in 1901 by Bouton \cite{Bou901},  
but only for a special graphs 
corresponding to the game of NIM. 
In 1944 this algorithm was extended to 
arbitrary digraphs by von Neumann and Morgenstern \cite{NM44}. 
In graph theory, the set of positions with even  $\R$  
is called the {\em kernel} of a digraph.
In the impartial game theory 
this set is referred to as the set {\em P-positions}, 
while the complementary set is the set of {\em N-positions}. 
These sets are characterized by the following two properties: 
There are no moves between P-positions and  
from every N-position there is a move to a P-position.

Smith's function $\R$ is a refinement 
of the concept of P-positions. 
It has the following, stronger, properties: 
\begin{itemize}   
\item[({\bf P})]
$x$ is a P-position if and only if  $\R(x)$  is even. 
In this case the player who makes a move from $x$ cannot win,  
yet, can resist for $\R(x)$ moves, but cannot guarantee a longer play. 
\item[({\bf N})]
Position $x$ is an N-position if and only if  $\R(x)$  is odd.  
Then the player making a move from  $x$ can win in at most 
$\R(x)$ moves, but cannot guarantee a faster win. 
\item[] 
In fact an optimal play is to   
reduce $\R$  by one, whenever it is positive.
\end{itemize}

\subsection{Sprague-Grundy (SG) functions}  
\label{ss01a}
Another refinement of the concept of P-positions 
was introduced by Sprague \cite{Spr36} and Grundy \cite{Gru39} 
as follows. 
Given an impartial game $\Gamma$, 
its SG function $\cG = \cG(\Gamma)$  is defined recursively: 
$$\cG(x) = mex\{\cG(y) \mid x \to y\},$$ 
where  $x \to y$  is a move from  $x$  to $y$  and 
mex  is the 
{\em minimum excludant}: 

\medskip 

$mex(S)$  is the minimum number in $\ZZ_+ \setminus S$
for every finite subset  $S \subset \ZZ_+$. 

\medskip 

In particular, $mex(\emptyset) = 0$ or, in other words, 
$\cG(x) = 0$  for any terminal position $x$  in $\Gamma$. 
The SG function is characterized by the following properties: 

\begin{itemize}   
\item[(I)] $\G(x) \neq \cG(y)$  for any move  $x \to y$.  
\item[(A)] For any position $x$  and integer $\ell$ 
such that $0 \leq \ell < \cG(x)$  
there exists a move  $x \to y$  such that 
$\cG(y) = \ell$.
\end{itemize} 

Furthermore, $x$  is a P-position if and only if  $\cG(x) = 0$. 

\subsection{Sprague-Grundy and Smith's Theories}  
\label{ss01}
Given $n$ impartial games 
$\Gamma_1, \dots, \Gamma_n$, 
by one move a player chooses one of them 
and makes a move in it. 
The player who has to move but cannot is the looser. 
The obtained game  
$\Gamma = \Gamma_1 \vee \dots \vee \Gamma_n$  
is called the {\em disjunctive compound} of games
$\Gamma_1, \dots, \Gamma_n$. For example, NIM with $n$ piles 
is the disjunctive compound of $n$  one-pile NIMs. 
The SG function $\G$ of $\Gamma$ is uniquely determined 
by the SG functions of $n$ compound games by formula 
$\G(\Gamma) = \G(\Gamma_1) \oplus \dots \oplus \G(\Gamma_n)$ 
where $\oplus$  is the so-called NIM-sum. 
These results were obtained by Bouton 
\cite{Bou901}  for the special case of NIM and 
then extended to arbitrary impartial games 
by Sprague \cite{Spr36} and Grundy \cite{Gru39}. 

\medskip 

Now suppose that by one move a player makes a move 
in each of $n$ compound games, rather than in one of them. 
Again, the player who has to move but cannot is the looser.
The obtained game  
$\Gamma = \Gamma_1 \wedge \dots \wedge \Gamma_n$  
is called the {\em conjunctive compound}. 
Remoteness function $\R$ of $\Gamma$ is uniquely determined 
by the remoteness functions of $n$ compound games by formula 
$\R(\Gamma) = \min(\R(\Gamma_1), \dots, \R(\Gamma_n))$. 
This result was obtained in 1966 by Smith \cite{Smi66}; 
see also \cite{Con76,Len07}.

\subsection{Main Results}

In Section \ref{s2} we consider Moore's NIM game \cite{Moo910}
and provide a simple, efficiently computable closed form expression for the remoteness of P-positions, and prove that computing the value for N-positions, is NP-hard, in general. 

In Section \ref{s2.5} we consider hypergraph NIM games (see \cite{BGHM15,BGHMM15}) and provide an efficiently computable closed form of the remoteness function for minimally transversal-free hypergraphs \cite{BGHMM17,BGHMM18}.

In Section \ref{s3} we consider the game Euclid introduced in \cite{CD69} and prove that the remoteness function can be computed in polynomial time for these games. 

In Section \ref{s4} we consider the game Wythoff \cite{Wyt907} and its generalizations \cite{Fra05,Gur12}. We provide an efficiently computable closed form for the remoteness value of P-positions, and a polynomial time algorithm for N-positions.  

\bigskip

\section{Moore's NIM}\label{s2}
This game was introduced in 1910 by Moore \cite{Moo910} as follows. 
Let $n$ and $k$ be two integers such that $0 < k \leq n$.   
Given $n$  piles of stones, two players alternate. 
By one move, a player chooses any $\ell$ non-empty piles with $1\leq \ell \leq k$ 
and reduces them arbitrarily (but strictly). 
If there exist no  $k$  non-empty piles, the game is over. 
Denote the obtained game by $\mim{n,k}$.
It was considered in several papers; 
see, for example, \cite{JM80,Niv04e,Niv04w,Niv04t} and also 
\cite{BGHM15,BGHMM15,BGHMM17,BGHMM18}. 

A position in $\mim{n,k}$ is 
$x=(x_1,\dots, x_n)\in\ZZ_+^n$, an $n$-dimensional 
vector with nonnegative integer entries.
Similarly to the standard NIM, 
we use the binary representations of integers $x_i$ 
\[
x_i = \sum_{j=0}^N 2^j x_{ij},
\]
We call the matrix of the $x_{ij}$ entries the \textit{Bouton matrix} of position $x$,
and define Moore's function $M(x)$  by
\[
M(x) = \sum_{j = 0}^N 
(k+1)^j \left[\sum_{i=1}^n x_{ij}\pmod{k+1}\right]. 
\]

Game $\mim{n,k}$ was solved by Moore as follows. 
\begin{theorem}[Moore \cite{Moo910}] \label{t-Moore}
A position $x$ is a \textup{P}-position 
in $\mim{n,k}$ if and only if  $M(x) = 0$. 
\end{theorem}
Note that by the above claim a position $x$ is a P-position if and only if 
\begin{equation}\label{e-Moore}
    \sum_{i=1}^n x_{ij} \equiv 0 \pmod{k+1} ~~ \mbox{ for all } j=0, \dots , N.
\end{equation}

\begin{remark}
In fact, the SG and Moore's functions are equal, 
$\G(x) = M(x)$, 
if one of them takes value at most 1; 
see \cite{JM80} and also \cite{BGHMM17} for more details.  
Yet, in general, no explicit formula 
or polynomial algorithm computing  the SG function is known for $\mim{n,k}$.
In particular, no closed formula is known already  
for the positions of the SG-value 2  for $n = 4$ and  $k =2$. 
\end{remark}

Let $S(x)$ denote the total number of stones 
in piles of a position $x$, that is, 
$S(x) = \sum_{i=1}^n x_i$.

\begin{theorem}\label{t-P-POSITIONS}
If $x$ is a \textup{P}-position in $\mim{n,k}$, then $\R(x) = 2S(x)/(k+1)$. 
\end{theorem}

\begin{proof} 
We claim first that in a P-position $x$ the sum $S(x)$ is a multiple of $k + 1$. 
Note that by \eqref{e-Moore} we have $\sum_{i=1}^nx_{ij}=\alpha_j(k+1)$ for some $\alpha_j\in\ZZ_+$ for all $j=0,...,N$. Thus we can write 
\[
S(x)=\sum_{i=1}^nx_i=\sum_{i=1}^n\sum_{j=0}^N2^jx_{ij}=\sum_{j=0}^N2^j\sum_{i=1}^n x_{ij}=\sum_{j=0}^N2^j\alpha_j(k+1),
\]
proving our claim.

Let us recall that the definition of $\R(x)$ assumes that both players play optimally. In particular, from a P-position the game returns to a P-position after two moves. 

Let us first prove that for every two consecutive moves the sum $S(x)$ is decreased by at least $k+1$. To see this, assume we move $x\to y$ and denote by $j$ the largest index such that $y_{ij}\neq x_{ij}$ for some $1\leq i\leq n$. Since we can decrease at most $k$ piles, we have $\sum_{i=1}^ny_{ij}=\sum_{i=1}^nx_{ij}-\delta$ for some $1\leq \delta \leq k$. 
Now our opponent must return to a P-position with a move $y\to z$. Let us denote by $\ell$ the maximal index such that $\sum_{i=1}^ny_{i\ell}\neq \sum_{i=1}^nz_{i\ell}$. Note that  $\ell\neq j$ implies that at least one of $\sum_{i=1}^nz_{i\ell}$ and  $\sum_{i=1}^nz_{ij}$ is not congruent to zero modulo $(k+1)$. 
Thus we have $\ell=j$. 
Furthermore, since $\sum_{i=1}^ny_{ij}\equiv (k+1)-\delta\pmod{k+1}$, the opponent must decrease at least $k+1-\delta$ piles to restore the condition that $\sum_{i=1}^nz_{id}\equiv 0\pmod{k+1}$ for all digits $d=0,\dots,N$. 
Thus in the first move at least $\delta$ piles were decreased, and in the response by the opponent at least $(k+1)-\delta$ piles  were changed. Thus we have $S(x)-S(z)\geq k+1$. 

We argue next that the player moving from a P-position has a move $x\to y$ that forces a response such that exactly $k+1$ stones are removed from the piles in these two moves.

Let $d$ denote the lowest digit such that $\sum_{i=1}^nx_{id}\neq 0$, and 
let $p$ be the smallest index of a pile such that $x_{pd}\neq 0$. We claim that moving to $y$, defined by $y_p=x_p-1$ and $y_i=x_i$ for all $i\neq p$, is such a move.
Namely, the opponent has a unique move $y\to z$ from $y$ to a P-position $z$ by \eqref{e-Moore}. Since $x$ is a P-position, we must have a subset $I\subseteq \{1,\dots ,n\}$ of cardinality exactly $k$ such that $y_{id}=1$ for all $i\in I$.  Defining $z_i=y_i-1$ if $i\in I$ and $z_i=y_i$ for $i\not\in I$ we get a P-position $z$ such that $y\to z$ is a move in which we remove exactly $k$ stones from the piles. Note that by \eqref{e-Moore}
all moves from $y$ to a P-position must have the same form for some index set $I$ as above.
\end{proof}

Let us next observe an important property of winning moves from N-positions.
For $x$ is an N-position $x$, and  a subset $K$ of the piles of cardinality $k$, we denote by $Y_K(x)$ the subset of P-positions that can be reached from $x$ by a move that decreases only some of the piles in $K$. 

\begin{lemma}\label{l-N-move}
Let $x$ be an \textup{N}-position, and 
$K \subseteq \{1,...,n\}$ be a subset of the piles of cardinality $k$.  
Then we have $S(y)=D$ for some constant $D$ and all positions $y\in Y_K(x)$. 
\end{lemma}

\begin{proof}
Since $K$ has cardinality less than $k+1$, 
for any two positions $y, y'\in Y_K(x)$ and any $j\leq N$ we must have $\sum_{i=1}^ny_{ij}=\sum_{i=1}^ny'_{ij}$ by \eqref{e-Moore},
from which the claim follows.
\end{proof}

For an N-position $x$ we have the following approach to compute $\R(x)$: 
find a winning move to a P-position of the player that minimizes 
the number of remaining stones; say, it reduces it to $D$. 
Then  $\R(x) = 1 + 2D/(k+1)$.

This yields a polynomial algorithm computing $\R(x)$ 
if $k$ or $n-k$ is fixed. 
Namely, we can consider all subsets $K$ of piles of cardinality $k$, 
and check if there is a winning move using only these piles. 
If there are many different winning moves using the same set $K$, 
then by Lemma \ref{l-N-move}, all of them results in the same number of stones. 
Thus by considering all $\binom{n}{k}$ subsets, which are polynomial if $k$ or $n-k$ is a constant,  and 
choosing the one for which the winning move 
leaves the minimum number of stones, we can obtain $D$ above.

We show next that computing the $\R$-function 
is a computationally hard task, in general.
More formally, consider the following decision problem \textup{R-NIM}: 
Given positive integers $n$, $k \leq n$, $x_1,\dots, x_n$, 
and $B$,  
we decide whether $\R(x)\leq B$ in $\mim{n,k}$. 

\begin{theorem}
 If $x$ is an \textup{N}-position in $\mim{n,k}$, then  \textup{R-NIM} is \textup{NP}-complete.
\end{theorem}

\begin{proof}
  Let $x$ be an N-position in $\mim{n,k}$ and 
  $a_j$ be the residue of $\sum_{i=1}^nx_{ij}$ modulo $k+1$. 
  If $a_j\ne 0$, a winning move of the first player 
  (a move to a P-position) should change it to 0 
  by changing elements of $j$-th column of the position matrix. 
  (Also, the change should obey an obvious restriction: 
  changing $0$ to $1$ is possible, if in some column 
  to the right of the current one $1$ is changed to $0$.)

  Suppose we have changed bits of $j$-th column 
  in rows forming a set $S$, $|S|\leq k$, and 
\[
u_j = \sum_{i\in S} x_{ij}
\]
be the number of ones in the corresponding cells of the matrix. 
If a move leads to a P-position, the number of ones in  $j$-th column 
should be a multiple of $k+1$. 
Thus, the number $u'_j$ of ones in $j$-th column and 
rows from $S$ after the move is given by the rule:
\[
u'_j = \left\{\begin{aligned}
&u_j-a_j, &&\text{if}\ u_j-a_j\geq0,\\
&k+1+u_j-a_j, &&\text{if}\ u_j-a_j<0.
\end{aligned}\right.
\]

The total number of stones is
\[
\sum_{ij}2^jx_{ij}.
\]
Therefore, the change $\Delta_j$ 
of the number of stones due to changes in $j$-th column is 
\begin{equation}\label{eq:Delta}
\Delta_j = \left\{\begin{aligned}
2^j&(-a_j), &&\text{if}\ u_j-a_j\geq0,\\
2^j&(k+1-a_j), &&\text{if}\ u_j-a_j<0.
\end{aligned}\right.
\end{equation}

Thus the number of stones after a move to a P-position is at least
\[
S(x) - \sum_{j}2^j a_j.
\]
 A move is called \emph{maximal} if it attains this bound. 

 Based on these observations, we introduce an auxiliary problem MAX-NIM: 
 Given positive integers $n$, $k \leq n$, and $x_1,\dots, x_n$ such that  $x=(x_1,\dots, x_n)$ is an N-position in $\mim{n,k}$, we decide whether
 there exists a maximal move from $x$. 

By Theorem \ref{t-P-POSITIONS}, a maximal move does exist if and only if
\[
\R(x) \leq 1+2\cdot\frac{S(x) - \sum_{j}2^j a_j }{k+1}\,.
\]
It gives a polynomial reduction of MAX-NIM to R-NIM.

To complete the proof, we need to show that MAX-NIM is NP-hard. For this purpose we reduce the problem Vertex Cover to MAX-NIM. 
An instance of Vertex Cover consists of a graph $G = (V,E)$ and an integer $c$. 
The question is whether $G$ contains  a vertex cover of size at most $c$. 

We associate a binary matrix $M_G$ to $G=(V,E)$ in the following way: Let $I_e$, $e\in E$ be pairwise disjoint sets of size $c$ such that they are also disjoint from $V$. The rows of $M_G$ are indexed by $V\cup \bigcup_{e\in E}I_e$, and the columns by $E$. For $v\in V$ and $e\in E$ we have a $1$ in the $(v,e)$ entry of $M_G$ if and only if $v\in e$. For a row $i\in I_f$, $f\in E$ and column $e\in E$ we have a $1$ in the $(i,e)$ entry of $M_G$ if and only if $f=e$. We define $n=|V|+c|E|$ and set $k=c$. We view the $n$ rows of $M_G$ as the binary encoding of $n$ positive integers $x_1, \dots , x_n$.

Let us observe that the sum of each columns of $M_G$ is $1$ modulo $k+1$, and thus by \eqref{e-Moore} the vector $x=(x_1,\dots,x_n)$ is an N-position. Thus, for this position by the above definitions we have $a_e=1$ for each columns $e\in  E$. 

We claim that this position and $k$ value is a YES instance of problem MAX-NIM if and only if the graph $G$ has a vertex cover of size at most $c$. 

Assume that $G$ has a vertex cover $C\subseteq V$ of size $|C|\leq c$. Then for each column $e\in E$ we can choose a $v\in C$ such that $v\in e$. For the corresponding entry $(v,e)$ of $M_G$ we switch the entry value $1$ to $0$. This reduces the column sums of $M_G$ by exactly one, and hence this change represents a maximal move from $x$ to a P-position. 

For the reverse, let us assume that there exists a maximal move to a P-position from $x$, and denote by $C$ the set of piles ($C\subseteq V\cup \bigcup_{e\in E}I_e$), in which we decrease the number of stones. Note that if there is a maximal move, then there is also one in which we switch exactly $a_e=1$ digit to $0$ in every column $e\in E$. Let us assume that  the considered maximal move has this property. 

Assume now that for some column $e=(u,v)$ this maximal move changes a bit in row $i\in I_e$. Note that in this row all other entries are $0$ by the definition of $M_G$. Let us now replace row $i$ by row $v$, if $v\not\in C$ and just delete row $i$ otherwise, and switch entry $(v,e)$ of $M_G$ to zero, instead of entry $(i,e)$. Since in this maximal move we have exactly one $1$ switched to $0$ in every column, the above change is feasible, and yields a new maximal move, with the same property. We can repeat this, until we get $C\subseteq V$, in which case it must be a set cover. Since in this process we do not increase the size of $C$, $C$ is a set cover of at most $c$ vertices. 

Since Vertex Cover is an NP-complete problem, the above proves that R-NIM is NP-hard, and since it belongs to NP, it is also NP-complete. 
\end{proof}

\medskip 

We can rephrase the above results as follows. 
It is easy to compute how to resists as long as possible from a P-position, 
but it is hard to compute how fast we can win from an N-position. 

\begin{remark}
Let us modify game  $\mim{n,k}$ as follows. 
By one move it is allowed 
to choose exactly  $k$  non-empty piles 
(rather than at most $k$)  and reduce 
each of them exactly by  $1$.  
The remoteness function of such game 
was efficiently computed 
in \cite{GMMV23} for $n=k+1$; see also \cite{GN23}. 
The mis\`ere version 
was considered in \cite{GMMN23}.
\end{remark}

\section{On Remoteness Functions of Hypergraph NIM}\label{s2.5}
Let us recall the so-called {\em hypergraph NIM games} \cite{BGHMM17}: 
Given a hypergraph $\cH \subseteq 2^{[n]} \setminus \{\emptyset\}$,  
the positions are integer vectors $x \in \ZZ_+^n$, 
and $x \to x'$ is a move if $x \geq x'$ and $\{i \mid x'_i < x_i\} = H\in \cH$. 
We call this move an {\em $H$-move} and denote it by $x\hmove{H}x'$. 
We call $x\hmove{H}x'$ a \textit{slow $H$-move}, 
if $x'_i= x_i-1$ for all $i\in H$.
Note that a position $x$ is terminal if and only if 
for every  $H \in \cH$  there exists an $i \in H$ such that $x_i = 0$. 

For instance, Moore's NIM is a hypergraph NIM corresponding to $\cH=\{H\subseteq [n]\mid 1\leq |H|\leq k\}$.

\medskip

Here we are concerned of determining the remoteness function 
for some special classes of hypergraph NIM games. 

Given a hypergraph $\cH\subseteq 2^{[n]}$, 
denote by $\cH^t$ the family of its transversals, i.e., 
$\cH^t=\{T\subseteq [n]\mid T \cap H \not=\emptyset  \mbox{ for all } H\in \cH\}$. 
For a subset $S\subseteq [n]$ we denote 
by $\cH_S$ the subhypergraph of $\cH$ induced by $S$, i.e., 

\[
\cH_S ~=~ \{H\in \cH \mid H\subseteq S\}.
\]

Let us consider a $NIM_\cH$ game for a hypergraph $\cH\subseteq 2^{[n]}$. 
To a position $x=(x_1,...,x_n)$ let us assign  
\[
\begin{array}{rl}
	m(x)&=\displaystyle \min_{i\in[n]} x_i, \text{ and }\\*[7mm]
	M(x)&=\{i\in [n]\mid x_i=m(x)\}.
\end{array}
\]

For an arbitrary nonnegative integer $k\in\ZZ_+$  define two sets of positions: 
\begin{eqnarray}
P(k) &=& \{ x\in\ZZ_+^n\mid m(x)=k \text{ and } M(x)\in\cH^{t} \}
\\
N(k) &=& \{ x\in\ZZ_+^n\mid ~\exists\; \text{a move } x\to x'\in P(k)\}. \label{eq-nk-mhf1} 
\end{eqnarray}

\begin{lemma}\label{l-P(k)}
	For all $k\in \ZZ_+$ and positions $x\in\ZZ_+^n$ with $m(x)\geq k$, we have $m(x')<k$ for all moves $x\to x'$ if and only if $x\in P(k)$.
\end{lemma}
\begin{proof}
	Consider an arbitrary position $x\in P(k)$, hyperedge $H\in\cH$, and $H$-move $x\hmove{H} x'$. Since $M(x)\in \cH^{t}$, we must have $H\cap M(x)\neq\emptyset$, and thus by the definition of an $H$-move we have $x'_j<x_j=m(x)$ for all $j\in H\cap M(x)$, implying $m(x')<m(x)$, as claimed.
	
	Let us next consider a position $x$ with $m(x)\geq k$ such that $x\not\in P(k)$. If $m(x)>k$, then we clearly have a move $x\to x'$ such that $m(x')\geq k$. If $m(x)=k$,  then we have $M(x)\not\in \cH^{t}$ by $x\not\in P(k)$. Thus there exists a hyperedge $H\in \cH$ with $H\cap M(x)=\emptyset$. Then for a slow $H$-move $x\to x'$ we get $m(x')=k$. 
\end{proof}

For a position $x\in\ZZ_+^n$ let us denote by $\cH(x)=\cH_{[n]\setminus M(x)}$.
Then $M(x)\not\in\cH^{t}$ implies 
\begin{equation}\label{e-1}
\{H\in\cH(x)\mid H\cup M(x)\in \cH^{t}\} ~=~ \cH(x)\cap \cH(x)^{t}.
\end{equation}
\begin{lemma}\label{l-N(k)}
If $\cH$ satisfies  $\cH\cap \cH^{t}=\emptyset$,  then we have
\[
N(k) ~=~ \left\{ x\in\ZZ_+^n\left| 
\begin{array}{c}
	m(x)=k,\\
	M(x)\not\in\cH^{t},\\
	\cH(x)\cap \cH(x)^{t} \neq \emptyset
\end{array} 
\right.\right\}. 
\]
\end{lemma}
\begin{proof}
	Let us arbitrarily take a position $x\in N(k)$ and an $H$-move $x\hmove{H} x'$ such that $x'\in P(k)$. 
	
Since  $x\to x'$ is a move,  we have $m(x')\leq m(x)$, and hence the relation $m(x)\geq k$ is implied by the definition of $N(k)$. 
We claim that $m(x)=k$. Assume for a contradiction that  $m(x)>k=m(x')$. Then we have $M(x')\subseteq H$. By $\cH\cap \cH^{t}=\emptyset$, we must have $M(x')\not\in\cH^{t}$, which is a contradiction with the assumption that $x'\in P(k)$. Thus we have $m(x)=k$. 
Furthermore, if $M(x)\in\cH^{t}$, then $M(x) \cap H\not=\emptyset$ holds by definition. This implies that 
$m(x')<k$, a contradiction. Hence we have $M(x)\not\in\cH^{t}$. 
To show that $\cH(x)\cap \cH(x)^{t}\neq \emptyset$, 
 we first note that  $H \cap M(x)=\emptyset$ and $M(x')\subseteq H\cup M(x)$, since $m(x)=m(x')=k$.  
By the definition of $P(k)$ we have $M(x')\in \cH^{t}$ and thus by \eqref{e-1} we have $H\in \cH(x)\cap \cH(x)^{t}$.

Let us finally consider an arbitrary position $x\in \ZZ_+^n$ such that $m(x)=k$, $M(x)\not\in \cH^{t}$, and $\cH(x)\cap \cH(x)^{t} \neq \emptyset$. Let us then choose a hyperedge $H\in \cH(x)\cap \cH(x)^{t}$ and consider the position $x^*\in\ZZ_+^n$ defined by $x^*_j=x_j$ for $j\notin H$ and $x^*_j=m(x)$ for $j\in H$. Now $x\to x^*$ is an $H$-move, and $M(x^*)=H\cup M(x)\in \cH^{t}$ by \eqref{e-1}. Thus we have $x^*\in P(k)$, implying that $x\in N(k)$. 
\end{proof}

Let us call a hypergraph \emph{minimally transversal-free} 
(MTF) if $\cH\cap \cH^{t}=\emptyset$ and 
for all nonempty proper subsets $S\subsetneq [n]$ there exists a hyperedge 
$H\in\cH_S$ such that $H\in(\cH_S)^{t}$.

There are several examples for MTF hypergraphs, 
including symmetric hypergraphs \cite{BGHMM15}, 
self-dual matroid hypergraphs \cite{BGHMM17}, 
exact-$k$-NIM games with $n=2k$ \cite{BGHMM15}, 
JM-games (e.g., Moore's game with $n = k+1$) \cite{BGHM15}, etc.

\bigskip

\begin{lemma}\label{l-N(k)-P(k)}
For an MTF hypergraph $\cH$,  we have
\[
N(k) ~=~ \left\{ x\in\ZZ_+^n \mid
	m(x)=k \text{ and }
	M(x)\not\in\cH^{t}
\right\}. 
\]
\end{lemma}
\begin{proof}
Note that condition $M(x)\not\in \cH^{t}$ implies $M(x)\neq [n]$ from which $\cH(x)\cap \cH(x)^{t} \neq \emptyset$ follows, since $\cH$ is minimally transversal-free. Thus the claim follows by Lemma \ref{l-N(k)}.
\end{proof}

\begin{corollary}\label{c-N(k)-P(k)}
Let $\cH\subseteq 2^{[n]}$ be an MTF hypergraph. Then for every integer $k\in \ZZ_+$, we have 
\[
N(k)\cup P(k) ~=~ \{x\in \ZZ_+^n\mid m(x)=k \}.
\]
\end{corollary}
\begin{proof}
The claim follows by the definition of $P(k)$ and Lemma \ref{l-N(k)-P(k)}.
\end{proof}

\begin{theorem}\label{t-S-MTF}
	Let $\cH\subseteq 2^{[n]}$ be an MTF hypergraph. Then we have the following relations for all $k\in \ZZ_+$: 
	\[
	\R(x) = 
	\begin{cases}
		2k &\text{ if } x\in P(k),\\
		2k+1 &\text{ if } x\in N(k)
	\end{cases}
	\]
\end{theorem}

\begin{proof}
We are going to prove this by induction on $k$. 

Note first that the terminals of the game 
are exactly the positions in $P(0)$, and 
$N(0)$ contains all positions from which a terminal can be reached. 
Thus the Smith values $0$ and $1$ are correctly assigned. 

Let us assume now that we already proved that the above rule describes the correct Smith value assignment for all $k<\ell$, that all positions in 
\[
\bigcup_{k=0}^{\ell-1} (P(k)\cup N(k))
\]
have their correct Smith valued assigned by the above rule. Let us remove all these positions and consider the residual game. Note that by Corollary \ref{c-N(k)-P(k)} we removed exactly the positions $x$ for which we have $m(x)<k$. Thus the residual game consists of the positions 
\[
\{x\in\ZZ_+^n\mid m(x)\geq k\}.
\]
Lemma \ref{l-P(k)} then implies that $P(k)$ 
is the set of terminals in this residual game, and thus by \eqref{eq-nk-mhf1}, 
$N(k)$ is exactly the set of positions from which $P(k)$ is reachable by a move. Consequently, the above formula assigns the correct Smith values to all positions in $P(k)\cup N(k)$.
\end{proof}

\begin{remark}
Moore's $\mim{n,k}$  is minimally transversal-free if and only if $n = k+1$.
In this case  for a P-position $x$ we have $\R(x) = 2S(x)/(k+1)$ by Theorem \ref{t-P-POSITIONS} which is in agreement to Theorem \ref{t-S-MTF}.
Indeed, for $x\in P(k)$ every column in the Bouton matrix is of size $n=k+1$  
and the sum is  $0 \pmod{k+1}$. 
Hence, all entries of any column are equal, 
implying $x_1 = \dots = x_n=m(x)=k$, since $x\in P(k)$. 
Thus, $2S(x)/(k+1)=2k$, as stated in Theorem \ref{t-S-MTF}.
\end{remark} 
\begin{remark}
Let us also add that for an MTF hypergraph and position $x\in\ZZ_+^n$ we can test in polynomial time if $x\in P(k)$ by the definition of $P(k)$ and membership $x\in N(k)$ is also testable in polynomial time by Lemma \ref{l-N(k)-P(k)}. 
Thus, by Corollary \ref{c-N(k)-P(k)},  
the value $\R(x)$ can be computed in polynomial time. 
\end{remark} 

\begin{remark}
A move $x \to x'$  can increase or keep the remoteness function 
value but only if  $x$ is an N-position. 
In other words, inequality $\R(x') \geq \R(x)$  
may hold only if $\R(x)$  is odd. 
\end{remark}

As a second example we consider the so called SG-decreasing games.

By definition of the SG function, 
no move  $x \to x'$  can keep its value, 
$\G(x) \neq \G(x')$.   
Furthermore, for any integer $m$ 
such that  $0 \leq m < \G(x)$ 
there exists a move  $x \to x'$  such that $\G(x') = m$. 

An impartial game is called {\em SG-decreasing} 
if $\G(x') < \G(x)$  for every move  $x \to x'$.  
These games were recently studied in \cite{BGHMM23}. 
In fact, a game is SG-decreasing if and only if
$\G(x) = h(x)$  for every position $x$, 
where $h(x)$  is the length of the longest play from  $x$, 
called the {\em height} of $x$, see \cite{BGHMM23}. 
By definition, a position is  terminal 
if and only if it is of height zero. 
Each such position can be reached by a move 
from any other position in an SG decreasing game.

\begin{proposition}
The remoteness function of an SG-decreasing game is defined by 
\[
\R(x)=
\begin{cases}
    1 & \text{ if } x \text{ is non-terminal},\\
    0 & \text{ otherwise}.
\end{cases}
\]
\end{proposition}
\begin{proof}
By the above properties P-positions are 
exactly the terminals of an SG-decreasing game, and 
from all other positions we can reach a terminal in one move.
\end{proof}
Thus SG-decreasing games provide a class 
in which the remoteness function is not increasing.

\medskip 

However, it is NP-complete  
to decide whether a hypergraph NIM is SG-decreasing 
and, even when it is,  it still remains NP-hard to compute 
the height $h(x) = \G(x)$ \cite{BGHMM23}.

\section{Game Euclid}\label{s3}
The SG function of this game is given 
by a simple explicit formula. 
Although  no such formula for 
the remoteness function is known, still it can be computed in polynomial time. 

Game Euclid was introduced in 1969 
by Cole and Davie \cite{CD69} and then 
considered in several papers; 
see 
\cite{Col05,Fra05,Gur07,Len03,Len04,Niv04e,Spi73} for example.

\medskip 

Positions of this game are all pairs  $(x,y)$  of positive integers. 
Two players alternate. 
If  $x=y$  then $(x,y)$  is a terminal position. 
Otherwise, by one move a player subtracts 
any multiple of the smaller number from the larger one  
such that the difference is still positive. 
Assume wlog that $x > y$,  
then,  $(x,y) \to (x - \ell y, y)$
is a move if and only if 
$\ell > 0$  and $x - \ell y > 0$. 

The classical Euclidean algorithm 
chooses the maximal such  $\ell$. 
Yet, for players this condition is waved. 
According to Euclid, the game terminates in 
a unique position $(z,z)$, where $z$ 
is the greatest common divisor of $x$ and $y$, denoted by  GCD$(x,y)$. 
 
 Without loss of generality, we can assume that $x$ and $y$ are co-prime, 
 that is, $z = 1$, since two games that begin 
 in $(x,y)$  and  $(mx,my)$  are equivalent, where $m$ denotes a positive integer.      
 The following properties are known: 
 
 \begin{itemize}
\item[(P)] P-positions are characterized by the inequalities 
$x < \phi y$  and  $y < \phi x$,    
where  $\phi = (1 + \sqrt{5})/2 = 1.618 \ldots$  
 is the golden ratio \cite{CD69}. 
\item[(SG)] The SG-function 
is determined in 2004 by Nivash \cite{Niv04e}: 
\[
\G(x,y) = \lfloor |x/y - y/x| \rfloor.
\]
\item[(PE)] There exists a unique move from a P-position \cite{Niv04e}. 
\item[(NE)]  
From an N-position,  there exists a unique move to a P-position \cite{Niv04e}. 
\item[(U)]  An optimal play from $(x,y)$ is unique. 
By definition, $\R(x,y)$  is the length of this play. 
\end{itemize} 

Obviously, (PE) and (NE) imply (U). This togehter with the following lemma shows that $\R(x,y)$ can be computed in   linear time. 

\begin{lemma}\label{lemma-euclid}
Let  $(x,y)$ be  a position  with  $x > y$.   
\begin{description}
\item[{\rm (i)}] If  $(x,y)$ is a \textup{P}-position, then 
the unique move from it to $(x',y)$ provides 
$\frac{x'}{x} \leq 1 - \frac{1}{\phi}$
\item[{\rm (ii)}] If  $(x,y)$ is an \textup{N}-position, then   
the unique move from it to a \textup{P}-position  $(x',y)$ provides 
$\frac{x'}{x} \leq \phi- 1$. 
\end{description}
\end{lemma}

\begin{proof}
For (i), we note that $x'=x-y$ by Property (P), and hence 
\[
\frac{x'}{x} \, < \, \frac{\phi y - y}{\phi y} \,= \,1-\frac{1}{\phi}. 
\]
For (ii),  $x$ can be represented by 
 $x=my+x'$, where $m$ is a positive integer and $\frac{y}{\phi} < x' < \phi y$. 
 Then  we have 
\[
\frac{x'}{x} \, = \, \frac{x'}{my+x'} \, < \, \frac{\phi y }{my + \phi y} \,\leq  \,\frac{\phi}{1 + \phi} \, = \,\frac{1}{\phi}. 
\]
\end{proof}

\begin{theorem} 
 $\R(x,y)$ can be computed in time linear in $\log (xy)$.  
\end{theorem} 

\begin{proof}
By Lemma \ref{lemma-euclid}
any move reduces  
the larger number by a constant factor, which implies the statement.  
\end{proof}

Yet, in contrast to the SG function,  no explicit formula is known 
for the remoteness function. 
Only positions of $\R$-value 0, 1, 2, or 3  have simple structure.  

\begin{proposition}
Let $x$ and $y$ be positive integers with GCD$(x,y) = 1$. Then, 
\begin{itemize}
\item[{\rm (0)}]
$\R(x,y)=0$  only for  $(x,y) = (1,1)$; 
\item[{\rm (1)}]
 $\R(x,y)=1$  if and only if either $x=1$ or $y=1$. 
 \item[{\rm (2)}]
$\R(x,y)=2$  if and only if $x>1, y>1$,   and  $|x-y|=1$.
 \item[{\rm (3)}]
 $\R(x,y)=3$  if and only if 
$x>1, y>1$, and 
$(x = my \pm 1$  or  $y = mx \pm 1)$ for some positive integer $m$. 
\end{itemize}
\end{proposition}

\proof 
(0): According to the rules of the game, 
position $(x,y)$ is terminal if and only if $x=y$. 
Moreover,  $x = y = 1$, since  GCD$(x,y) = 1$. 
Thus, $\R(x,y) = 0$  if and only if  $x = y = 1$. 

(1): By definition, $\R(x,y) = 1$  if and only if there is 
a move  $(x,y) \to (1,1)$. 
Obviously, it exists if and only if 
$x=1$ or $y=1$, but not both. 

(2): By definition,  $\R(x,y) = 2$ 
if and only if 
there exists a move $(x,y) \to (x',y')$. 
and  $\R(x',y') = 1$ for each such move.    
Obviously, this happens if and only if 
$x > 0, y > 0$, and  $y = x \pm 1$. 

(3): By definition, $\R(x,y)=3$  if and only if 
there exists a move from $(x,y)$  
to a position of $\R$-value 2, 
and there is no move to a position of $\R$-value 0, 
that is, to the terminal position $(1,1)$. 
This happens if and only if 
$x>1, y>1$, and 
($x = m y \pm1$  or $y = mx \pm 1$ for some $m > 0$). 
\qed 

\medskip 

In general, $\R$-values look chaotic. 
For example, we have 
\[
\R(17,11) = 4, \; \R(17,12) = 6, \; \R(17,13) = 4,   
\]
and the corresponding optimal plays are:  
\begin{eqnarray*}
&&(17, 11) \to (6, 11) \to (6, 5 ) \to (1, 5) \to (1, 1); \\
&&(17, 12) \to (5, 12) \to (5, 7) \to (5, 2) \to 
(3, 2) \to (1,2) \to (1,1);\\ 
&&(17, 13) \to (4, 13) \to (4, 5) \to (4, 1) \to (1, 1).
\end{eqnarray*}

\section{Game Wythoff and Its Generalizations}\label{s4}
Consider the following version 
of NIM with two piles. 
Positions are pairs of nonnegsative numbers  $(x,y)$.  
Fix two positive parameters $a$ and $b$.  
Two players alternate. 
By one move it is allowed 
to reduce $x$ by $\epsilon$ and $y$  by  $\delta$ 
such that  $0 \leq \epsilon \leq x$, $0 \leq \delta \leq y$, 
and $\epsilon+\delta > 0$. 
Furthermore, we require that  either
$|\epsilon - \delta| < a$  or  $\min(\epsilon,\delta) < b$, or both.  

This game  WYT$(a,b)$  was introduced in \cite{Gur12}, 
case $b=1$  was considered in 1982 by Fraenkel \cite{Fra82},  
and case $a=b=1$, as early as in 1907 by Wythoff \cite{Wyt907}. 
There are several works related to Wythoff's NIM  
and its different generalizations; see for example
\cite{BF90,Con76,Cox53,FB73,Fra82,Fra84,
Fra05,Gur12,Niv04w,Niv04t,Wyt907,YY67}. 

Explicit formulas for the P-positions  
were obtained in \cite{Fra82,Wyt907}. 
Interestingly, no explicit formula or polynomial 
algorithm is known for computing the SG function, 
even when $a=b=1$, 
in spite of quite intensive research, 
which is summarized in \cite{Niv04w}.

In this section, we obtain a polynomial algorithm 
computing the remoteness function for WYT$(a,b)$, 
for any positive constant parameters $a$ and $b$.

Let $\{(x_m,y_m)\mid m \in \ZZ+\}$ 
be the set of P-positions of WYT$(a,b)$ such that $x_m \leq y_m$
for all nonnegative integers $k$.
Assume that $x_0 \leq x_1 \leq \dots$ holds. 
Then we 
note that $(x_0,y_0)=(0,0)$ and $x_m < x_{m+1}$ holds 
for all nonnegative integers $m$, 
by the rule of game WYT$(a,b)$. Furthermore, 
the following recursion  was obtained in \cite{Gur12}:  
\[
x_m = mex_b\{x_0, y_0, \dots, x_{m-1}, y_{m-1}\}  
\mbox{ and } y_m = x_m + am  \ \ \mbox{ for } m \in\ZZ_+.
\]
For a  finite set of nonnegative integers 
 $S = \{s_1, \dots, s_\ell\}$  with 
$s_1 < \dots < s_\ell$, 
the minimum $b$-excludant $mex_b(S)$ is defined by 
\[
mex_b(S)=
\begin{cases}
    0 & \text{if } S=\emptyset ,\\
    s_i + b & \text{if $i$  is the smallest index with 
$s_{i+1} - s_i > b$}, 
\end{cases}
\]
where we regard  $s_{\ell+1}=+\infty$. 
Clearly, $mex_1$  is the standard minimum excludant $mex$, \cite{Fra82}. 

By symmetry, 
 $(x,y)$ is a P-position if and only if so is $(y,x)$. 
Examples of P-positions can be found  
in \cite{Fra82} for $b=1$  and 
in \cite{Gur12} for  $a > 1$ and $b > 1$. 
For $a=b=1$ the above recursion was solved by Wythoff 
who proved that  
$x_m = \lfloor \phi m \rfloor$  and $y_m = x_m + m$. 
For $b=1$ it was solved by Fraenkel  
who proved that  
$x_m = \lfloor \frac{m}{2} \big(2-a + \sqrt{a^2 + 4}\big) \rfloor$ 
and $y_m = x_m + am$. 
For example, $x_m = \lfloor \phi m \rfloor$  for $a = 1$ and 
$x_m = \lfloor \sqrt{2} m\rfloor$  for $a = 2$.
Let us also recall from \cite{Gur12} that the following inequalities hold for any index $m\geq 0$:
\begin{equation}\label{e-Gur12}
b\leq x_{m+1}-x_m\leq 2b ~~\text{ and }~~ b+a\leq y_{m+1}-y_m\leq 2b+a.
\end{equation}

For positive constants $a$  and $b$, 
a polynomial algorithm computing P-positions 
(based on the Perron-Frobenius theory) 
was constructed in \cite{BGO13}, 
although no explicit formula for P-position was given.
More precisely, for given positive integers $M$, $X$ and $Y$, we can solve in polynomial time the following problems:
\begin{itemize}
\item[(i)] Determine $x_M$.
\item[(ii)] Find the maximum integer $m$ such that $x_m\leq X$.
\item[(iii)] Find the maximum integer $m$ such that $y_m\leq Y$. 
\end{itemize}

\medskip
In this game, we consider the following types of moves: 

\begin{description}
\item[{\rm (a)}] 
the numbers of taken stones can differ by at most $a-1$.
\item[{\rm (b)}]  reduce one pile arbitrarily, but  
the other one by at most $b-1$ stones. 
\end{description}
Moves of type (a) are called {\em diagonal}, 
while moves of type (b) are called 
{\em horizontal} (respectively, {\em vertical})  
if $x$ (respectively, $y$) 
is reduced by at most  $b-1$. 
Note that the same move may be of two 
or even three distinct types.

\begin{theorem}\label{t-WYT-P}
The remoteness function 
of a \textup{P}-position $(x_m,y_m)$  satisfies 
$\R(x_m,y_m) = 2m$  for all $m \in\ZZ_+$.
\end{theorem}

\begin{proof}
Consider the vertical move $(x_m,y_m) \to (x_m,y_m-1)$. 
By the characterization of P-positions, 
this cannot be a P-position, and hence it is an N-position. 

It is easy to see that from $(x_m,y_m-1)$ we can reach $(x_{m-1},y_{m-1})$ by one diagonal move.
We claim that from $(x_m,y_m-1)$ we cannot reach by any other P-position $(x_\ell,y_\ell)$ or $(y_\ell,x_\ell)$ with $\ell \not=m-1$, which implies that $\R(x_m,y_m-1)=\R(x_{m-1},y_{m-1})+1$.
Since $\R(x_{0},y_{0})=0$, this completes the proof. 

There are two cases to consider, if the claim is not true.

Assume for a contradiction first that there is a move $(x_m,y_m-1)\to (x_\ell,y_\ell)$ for some $\ell\not=m-1$. 
Since $(x_\ell,y_\ell)< (x_m,y_m-1)$, we have $\ell < m-1$ by \eqref{e-Gur12}. 
The changes in the coordinates are given by 
\[
\begin{array}{rl}
\epsilon & = x_m-x_\ell\\
\delta & = y_m-1-y_\ell \\
& = x_m-x_\ell +(m-\ell)a-1.
\end{array}
\]
Thus, we have $\epsilon <\delta$ for any $a\geq 1$. For such a move to be legal in WYT$(a,b)$, we need either $|\epsilon-\delta|<a$ or $\min\{\epsilon,\delta\} <b$. Here we have 
\[
|\epsilon-\delta|=(m-\ell)a-1 \geq 2a-1 \geq a
\]
for all $a\geq 1$. Furthermore, we have
\[
\min\{\epsilon,\delta\} =\epsilon =x_m-x_\ell\geq b
\]
by \eqref{e-Gur12}. Thus no such move exists.

Assume for a contradiction next that there is a move $(x_m,y_m-1)\to (y_\ell,x_\ell)$ for some $\ell \not=m-1$. 
Since $y_\ell\leq x_m$, we have $\ell < m-1$ by \eqref{e-Gur12}. 

The changes in the coordinates are in this case
\[
\begin{array}{rclcl}
\epsilon & =& x_m-y_\ell&=&x_m-x_\ell-\ell a\\
\delta & =& y_m-1-x_\ell&=&x_m-x_\ell +ma-1.\\
\end{array}
\]
Thus we have $\epsilon<\delta$ for any $a\geq 1$. For such a move to be legal in WYT$(a,b)$, we again need either $|\epsilon-\delta|<a$ or $\min\{\epsilon,\delta\} <b$. Here we have 
\[
|\epsilon-\delta|=(m+\ell)a-1\geq  a
\]
for all $a\geq 1$. Furthermore, $\min\{\epsilon,\delta\}=\epsilon=x_m-y_\ell$. By the definition of $x_m=mex_b\{x_i,y_i\mid 0\leq i<m\}$, we must have 
\[
\min\{\epsilon,\delta\}=\epsilon=x_m-y_\ell \geq b,
\]
which proves that no such move exists.  
\end{proof}

For N-positions we do not have a closed formula for the remoteness value.

We show in the rest of this section that the $\R(x,y)$ values can be computed in polynomial time for all positions $(x,y)$ of the WYT$(a,b)$ game. 

First we prove that from an N-position one can reach at most six P-positions by a move. Furthermore, these at most six positions can be determined in polynomial time.

\begin{lemma}\label{l-diagonal}
For an \textup{N}-position $(x,y)$ of WYT$(a,b)$, 
there are at most two P-positions 
that may be reachable from $(x,y)$ by a diagonal move. 
\end{lemma}

\begin{proof}
Assume w.l.o.g. that $x\leq y$ and write $y-x=a\ell +r$, for some $0\leq r < a$. 
Assume also that $(x,y)\to (x-\epsilon,y-\delta)$ is a diagonal move to a P-position. 
In particular, $\epsilon,\delta\geq 0$, $|\epsilon-\delta|<a$, and $\epsilon+\delta>0$. 
By the characterization of P-positions, we have $(y-\delta -(x-\epsilon) = 0 \pmod{a}$. Consequently, we have $\epsilon-\delta\in \{-r,a-r\}$.
Therefore, we have $x-\epsilon\leq y-\delta$,  and hence 
the only two candidate P-positions are $(x_\ell,y_\ell)$ or $(x_{\ell+1},y_{\ell+1})$. 
\end{proof}

Let us remark that since the above proof computes the index $\ell$ in polynomial time, the P-positions $(x_\ell,y_\ell)$ and $(x_{\ell+1},y_{\ell+1})$ are also computable in polynomial time by the above cited result (i) of \cite{BGO13}. 

\begin{lemma}\label{l-non-diagonal}
Given an \textup{N}-position $(x,y)$ of WYT$(a,b)$, 
there are at most four \textup{P}-positions that may be reachable 
from $(x,y)$ by a non-diagonal move. 
\end{lemma}

\begin{proof}
Assume w.l.o.g. that $x\leq y$ and  $(x,y)\to (x_\ell,y_\ell)$ is a vertical move to a P-position, where $x-x_\ell<b$. We claim that no other $(x_m,y_m)$ is reachable from $(x,y)$, since otherwise $|x_\ell-x_m|<b$ would follow, contradicting 
\eqref{e-Gur12}. Assume next that $(y_m,x_m)$ is also reachable from $(x,y)$ by a vertical move. 
Note that  $(y_m,x_m)$ is different from  $(x_m,y_m)$.  
Analogously to the previous case, it is implied that no other $(y_s,x_s)$ is reachable from $(x,y)$ by another vertical move. 

We can similarly prove that there are at most two P-positions that are reachable from $(x,y)$ by horizontal moves. 
\end{proof}

Let us remark that by using the above cited algorithmic results (ii) and (iii) of \cite{BGO13}, we can compute the largest $\ell$ such that $x_\ell<x$ or $y_\ell<x$ and if they are at most $b-1$ smaller than $x$, then using (i) we can compute the corresponding P-positions and check if those are indeed reachable form $(x,y)$, all in polynomial time. We can check analogously in polynomial time the P-positions reachable by horizontal moves.

\begin{theorem}\label{t-WYT-complexity}
    For a position $(x,y)$ of WYT$(a,b)$ we can compute $\R(x,y)$ in time polynomial in $a$, $b$ and $\log(xy)$. 
\end{theorem}

\begin{proof}
    Assume w.l.o.g. that $x\leq y$ and consider $\Delta=(y-x) \pmod{a}$.

    If $\Delta \neq 0$ then $(x,y)$ is an N-position. 

    If $\Delta=0$, then denote by $m=\frac{y-x}{a}$ and compute $x_m$ by the above cited polynomial algorithm (i). If $x=x_m$ then $(x,y)=(x_m,y_m)$ and $\R(x,y)=2m$ by Theorem \ref{t-WYT-P}. Otherwise $(x,y)$ is an N-position.

    From an N-position $(x,y)$ we can reach at most six P-positions, and by Lemmas \ref{l-diagonal} and \ref{l-non-diagonal} and the following remarks we can determine those P-positions in polynomial time. Then $\R(x,y)$ is one larger than the smallest of the corresponding remoteness values by the definition of $\R$.
\end{proof}

{\bf Acknowledgements.}
Research of the second and fourth authors was 
prepared within the framework 
of the HSE University Basic Research Program.  
The third author was partially supported by JST ERATO Grant Number JPMJER2301 and JSPS KAKENHI Grant Numbers JP20H05967, JP19K22841, and
JP20H00609.

\end{document}